\documentclass[leqno,a4paper,12pt]{article}
\usepackage{amsfonts,amssymb}

\font\ecal=eusm10 at 12pt
\def\esm#1{\hbox{\ecal {#1}}}

\def\Z{\mathbb{Z}}
\def\Q{\mathbb{Q}}

\def\P{\mathbb{P}}

\def\O{{\cal O}}

\def\til#1{\widetilde{#1}}
\def\ovl#1{\overline{#1}}
\def\pf{{\indent\textit{Proof.}\ }}
\def\qed{\hfill$\square$}

\def\PGL{\mathrm{PGL}}
\def\tree{\esm{T}}
\def\limitset{\esm{L}}
\def\fixset{\esm{F}}
\def\valuation{\nu}

\def\Vert{\mathrm{Vert}}
\def\Edge{\mathrm{Edge}}
\def\Ends{\mathrm{Ends}}
\def\Star{\mathrm{Star}}
\def\longhookrightarrow{\lhook\joinrel\longrightarrow}

\def\bigdownarrow{\vphantom{\bigg|}\Big\downarrow}
\newcounter{para}[section]
\renewcommand{\thepara}{\thesection.\arabic{para}}
\renewcommand{\thesection}{\arabic{section}}
\renewcommand{\paragraph}{\refstepcounter{para}
\indent{\bf\thepara.}}
\newcommand{\sectioning}{\refstepcounter{section}
\indent{\bf \thesection.}}

\newenvironment{condlist}{\vspace{1ex}
     \begin{list}{}{%
          \setlength{\topsep}{0pt}%
          \setlength{\parsep}{0pt}%
          \setlength{\itemsep}{0pt}}%
     }{\end{list}\vspace{1ex}}

\def\msq23{\sqrt{6}}
\def\m2sq3{2\sqrt{3}}

\title{{\large {\bf Graph Theoretic Construction of Discrete Groups\\
over $p$-adic Fields}}}
\author{Fumiharu Kato}
\date{}
\begin{document}
\maketitle
%

\sectioning\label{section-introduction}
{\bf Introduction}

\vspace{1ex}
Let $K$ be a finite extension of $\Q_p$ and $\Gamma$ a finitely
generated discrete subgroup in $\PGL(2,K)$. 
It is well-known (e.g.\ \cite[\S 1]{Her80}) that such a group $\Gamma$
contains a free normal subgroup of finite index (so called {\it Schottky
group}). Hence it gives, through Mumford uniformization
of an analytic curve, a Galois covering, possibly with ramifications,
over a 
projective curve $X_{\Gamma}$ with the Galois group $\Gamma$.
Like that Mumford uniformization is linked with the corresponding
geometry of (a subtree of) Bruhat-Tits tree, the situation as above is
well described by the action of $\Gamma$ on such trees (cf.\
\cite{vdP97}, \cite{CKK99}, and Proposition \ref{pro-branch} below).
For instance, the number of branch points and the branching degrees can
be calculated only by looking at the corresponding trees.
In slightly more precise terms, the discrete group $\Gamma$ gives rise to 
a certain subtree $\tree^{\ast}_{\Gamma}$ in the Bruhat-Tits tree acted
on by $\Gamma$, and the quotient graph
$T^{\ast}_{\Gamma}=\Gamma\backslash\tree^{\ast}_{\Gamma}$ can be
decorated to a {\it graph of groups}
$(T^{\ast}_{\Gamma},\Gamma_{\bullet})$, that is, a graph on which finite 
groups are attached to vertices and edges in a compatible way. 
There are several nice aspects in it; for instance, the genus of
$X_{\Gamma}$ is the first Betti number of $T^{\ast}_{\Gamma}$, ends of 
$T^{\ast}_{\Gamma}$ are in bijection with branch points which
preserves the decomposition groups, $\Gamma$ is isomorphic to the
direct-limit group (essentially by amalgams) associated to it, etc.

In this paper we will discuss the converse of the above procedure in
case $X_{\Gamma}$ is a rational curve. More precisely, we will answer to 
the following 

\vspace{1ex}
{\bf Question.} {\sl Given an abstract tree of groups
$(T,G_{\bullet})$, when can one find a finitely generated discrete
subgroup $\Gamma$ in $\PGL(2,K)$ over some $K$, isomorphic to the
amalgam group associated to $(T,G_{\bullet})$, such that the
tree of groups $(T^{\ast}_{\Gamma},\Gamma_{\bullet})$
``essentially'' coincides with the original $(T,G_{\bullet})$, i.e., 
roughly speaking, when is 
$(T,G_{\bullet})$ realizable in the above context?}

\vspace{1ex}
The meaning of ``essentially'' is that these trees of groups are the
same modulo finite subtrees which are, so to speak, futile parts both 
topologically and group theoretically; this can be precisely stated by
the notion of {\it contraction} (cf.\ \cite[Prop.\ 1]{CKK99} and
Definition \ref{def-contraction} below).

To answer the question, we will introduce the so-called {\it
$\ast$-admissibility} for 
such abstract trees of groups; roughly speaking, a tree of groups is
$\ast$-admissible if and only if it has nice compatible embeddings into
the Bruhat-Tits tree and $\PGL(2,K)$ satisfying a certain local
condition, local at each vertex. What we will prove is that this
actually gives the necessary and sufficient condition for the
realizability. 

Moreover, this condition leads to a practical way of constructing
discrete subgroups. For instance, suppose we are interested in
classifying all the possible such groups with given number of branch
points and braching degrees. Our method basically reduces the problem 
into a combinatorial business which is often easy in principle. 
In \S\ref{section-application}, we will exhibit two examples of such
constructions. Our method actually has more applications (e.g.\ $p$-adic
analogue of triangle groups, quadrangle groups, etc.), some of which will
be discussed elsewhere (cf.\ \cite{Kat00}).

It should be noted that, as one finds in the first example in
\S\ref{section-application}, our method of construction may be
viewed as a paraphrase of the known method by isometric circles (e.g.\  
\cite{Her80}). But the advantage of ours is its clear link with
the tree which enable us to construct more complicated groups.
In the second example in 
\S\ref{section-application} we will construct a series of diadic
triangle groups, which actually answers affirmatively to Yves Andr\'e's
expectation that there should be infinitely many $p$-adic {\it
non-arithmetic} triangle groups (cf.\ \cite[9.4]{And98}). 

The plan of this paper is as follows: In the next section
(\S\ref{section-tree}) we will collect basic notions such as Bruhat-Tits 
tree, subtrees associated to discrete groups, and trees of groups,
etc. (This section contains nothing new.) The $\ast$-admissibility will
be introduced in \S\ref{section-realization}, where our main theorem
will be proved. The final section \S\ref{section-application} exhibit
applications.

\vspace{1ex}
{\sl Notation and conventions}.\ 
Throughout this paper $K$ denotes a finite extension of $\Q_p$,
$\O_K$ the integer ring, and $\pi$ a prime element in $\O_K$. 
We write $[K\colon\Q_p]=ef$, where $e$ is the ramification degree and
$q=p^f$ is the the number of elements in the residue field
$k=\O_K/\pi\O_K$.
We denote by $\valuation\colon K^{\times}\rightarrow\Z$ 
the normalized (i.e., $\valuation(\pi)=1$) valuation. 

\vspace{2ex}
\sectioning\label{section-tree}
{\bf Trees and groups}

\vspace{1ex}
This section contains nothing essentially new. The notions and
statements in this section can be found in the references listed in the
end of this paper.
Proofs are put for the reader's convenience, and are often sketchy. 

\vspace{1ex}
\paragraph\label{para-bttree}\
{\sl Bruhat-Tits tree.}\ 
First we recall the basic properties of Bruhat-Tits tree $\tree_K$
attached to $\PGL(2,K)$. It is the tree whose vertices are similarity
classes of $\O_K$-lattices in $K^2$, and two vertices are connected by
an edge if the corresponding quotient module has length one. 
There is a canonical action by $\PGL(2,K)$ on $\tree_K$.
For a vertex $v$, which is the similarity class of $M\subset K^2$, 
edges emanating from $v$ are in canonical bijection with the lines
in $M/\pi M\cong k^2$, i.e., $k$-rational points of
$\mathrm{Proj}\,\mathrm{Sym}_k(M/\pi M)\cong\P^1_k$:
$$
\left\{
\begin{minipage}{5.7cm}
\setlength{\baselineskip}{.85\baselineskip}
\begin{small}
Edges $\sigma$ in $\tree_K$ emanating from 
$v$
\end{small}
\end{minipage}
\right\}
\longleftrightarrow
\P^1(k).
\leqno{\indent\textrm{{\rm (\ref{para-bttree}.1)}}}
$$
The set of ends (i.e.\ equivalence classes of half-lines different by a
finite segment) are canonically idetified with $K$-rational points of
$\P^1_K$, since they are ``limits'' of sequences of lattices with
length one successive quotients.
$$
\left\{
\begin{minipage}{2.2cm}
\setlength{\baselineskip}{.85\baselineskip}
\begin{small}
 Ends in $\tree_K$
\end{small}
\end{minipage}
\right\}
\longleftrightarrow
\P^1(K).
\leqno{\indent\textrm{{\rm (\ref{para-bttree}.2)}}}
$$
Note that this bijection is equivariant with the action by $\PGL(2,K)$.

\vspace{1ex}
\paragraph\label{para-treenotation}\ 
{\bf Notation.}\ 
For an abstract tree $T$ we denote by $\Vert(T)$ (resp.\
$\Edge(T)$, $\Ends(T)$) the set of all vertices (resp.\
unoriented edges, ends).
The notation $v\vdash\sigma$ for $v\in\Vert(T)$ and
$\sigma\in\Edge(T)$ means that $\sigma$ emanates from $v$. 
For a vertex $v\in\Vert(T)$ we denote by $\Star_v(T)$ the set of 
edges in $\Edge(T)$ emenating from $v$.
For two vertices $v_0$ and $v_1$, we denote by $[v_0,v_1]$ the geodesic
path connecting them.
For $\varepsilon_0,\varepsilon_1\in\Ends(T)$ and $v\in\Vert(T)$,
the unique straight-line (resp.\ half-line) connecting $\varepsilon_0$
and $\varepsilon_1$ (resp.\ $v$ and $\varepsilon_0$) is denoted by
$]\varepsilon_0,\varepsilon_1[$ (resp.\ $[v,\varepsilon_0[$).
The geometric realization $|T|$ is metrized so that the path
$[v_0,v_1]$ ($v_0,v_1\in\Vert(T)$) is of length equal to the number
of edges in it. The metric function is denoted by
$d_T(\cdot,\cdot)$, or simply by $d(\cdot,\cdot)$.
If $T\subseteq\tree_K$, then we always regard the set
$\Ends(T)$ as a subset of $\P^1(K)$ by (\ref{para-bttree}.2).

\vspace{1ex}
\paragraph\label{lem-compact}\
{\bf Lemma.}\ {\sl Let $\tree$ be a subtree of $\tree_K$. Then the set
of ends of $\tree$, regarded as a subset in $\P^1(K)$, is a closed
(hence compact) set.} 

\vspace{1ex}
\pf
Let $\{\varepsilon_n\}_{n=1}^{\infty}$ be a set of ends of $\tree$ which 
converges, as points in $\P^1(K)$, to a point $\varepsilon$.
What to prove is that $\varepsilon$ is contained in $\Ends(\tree)$.
For each $n$, let $u_n$ be the vertex of $\tree$ determined by 
$]\varepsilon_{n-1},\varepsilon_n[\,\bigcap\,]\varepsilon_n,
\varepsilon_{n+1}[\,=\,]\varepsilon_n,u_n]$.
Then the union of all the segments $[u_{n-1},u_n]$ in $\tree$ contains a
half-line $\ell$ pointing to the end $\varepsilon$.
\qed

\vspace{1ex}
\paragraph\label{para-compact}\
{\sl Tree from a compact set.}\ 
Next we recall the definition of trees from compact sets
(\cite[(2.4)]{CKK99}): 
Let $\limitset$ be a compact subset of $\P^{1,\mathrm{an}}_K$.
We assume that every point in $\limitset$ is at most $K$-valued.
The {\it tree generated by $\limitset$}, denoted by $\tree(\limitset)$,
is the minimal subtree in $\tree_K$ having $\limitset$ as the set of
ends; it is an empty tree if $\limitset$ consists of less than $2$
points.
This notion depends on the base field $K$, but differs only by subdivision.
Note also that the tree $\tree(\limitset)$ in general differs from the one
by Gerritzen-van der Put \cite[I.\S 2]{GvP80}; for instance, the tree
$\tree^{\mathrm{GvdP}}(\limitset)$ by them is a finite tree for $\limitset$
being finite, whereas ours are not. In fact, we have the following
criterion: 

\vspace{1ex}
(\ref{para-compact}.1) {\sl The tree $\tree^{\mathrm{GvdP}}(\limitset)$
coincides with $\tree(\limitset)$ if and only if
$\Ends(\tree^{\mathrm{GvdP}}(\limitset))=\limitset$.}

\vspace{1ex}\noindent
This can be easily seen by the fact that $\tree(\limitset)$ is the
minimal subtree containing all the apartments $]z,w[$ for
$z,w\in\limitset$ ($z\neq w$).

\vspace{1ex}
\paragraph\label{para-recall}\ 
{\sl Elements in a discrete subgroup.}\ 
The following facts are well-known, but are inserted herein for the
reader's convenience:
An element $\gamma\in\PGL(2,K)$ is said to be {\it parabolic} (resp.\
{\it elliptic}, resp.\ {\it hyperbolic}) if it has only one eigenvalue 
(resp.\ two distinct eigenvalues with equal valuations, resp.\ two
distinct eigenvalues with different valuations).
Let $\Gamma$ be a discrete subgroup of $\PGL(2,K)$.
Then:

\vspace{1ex}
(\ref{para-recall}.1) {\sl There exists no parabolic element in $\Gamma$ 
other than $1$.}

\vspace{1ex}
(\ref{para-recall}.2) {\sl An element $\gamma\in\Gamma$ is of finite 
order if and only if it is elliptic.}

\vspace{1ex}
Suppose two elements $\theta$ and $\chi$ have exactly one common fixed
point $\infty\in\P^1(K)$.
We may assume $\theta={u\ \,0\ \ \,\choose \,0\ \,u^{-1}}$ and $\chi={a\
\,b\ \ \,\choose \,0\ \,a^{-1}}$ with $b\neq 0$. Then 
it is easy to see that $\theta\chi\theta^{-1}\chi^{-1}$ is a parabolic
element. Hence: 

\vspace{1ex}
(\ref{para-recall}.3) {\sl No two elements in $\Gamma$ have exactly one 
common fixed point in $\P^1_K$.}

\vspace{1ex}
\paragraph\label{lem-discrete}\
{\bf Lemma.}\ {\sl
Let $\Gamma$ be a subgroup in $\PGL(2,K)$ acting on a subtree
$\tree$ of $\tree_K$.

(1) If $\Gamma$ is discrete, then for each $v\in\tree$ the stabilizer
of $v$ is a finite group.

(2) Conversely, if the stabilizer of at least one vertex $v$ is finite, 
then $\Gamma$ is discrete.}

\vspace{1ex}
\pf
(1) is well-known (the stabilizer in $\PGL(2,K)$ of a vertex is an open
compact subgroup.)
Suppose that there is a sequence $\{\gamma_i\}\subset\Gamma$ converging to
$1$, then, except for finitely many $\gamma_i$'s, they are contained in the
stabilizer of $v$, since the stabilizer of a vertex in $\PGL(2,K)$
is an open neighborhood of $1$.
\qed

\vspace{1ex}
\paragraph\label{para-discretegroup}\
{\sl Trees from a discrete group.}\ 
Let $\Gamma$ be a finitely generated discrete subgroup in $\PGL(2,K)$.
We may assume, replacing $K$ by a finite extension if necessary, that
every element ($\neq 1$) in $\Gamma$ has at most $K$-valued fixed points 
in $\P^1_K$ (cf.\ \cite[I.3.1\ (1)]{GvP80}).
Then $\Gamma$ acts on $\tree_K$ without inversion.
Let
\begin{eqnarray*}
\limitset_{\Gamma}&=&\textrm{the set of limit points of $\Gamma$},\\
\fixset_{\Gamma}&=&\textrm{the set of fixed points of elements 
($\neq 1$) in $\Gamma$}.
\end{eqnarray*}
These are subsets in $\P^{1,\mathrm{an}}_K$ satisfying
$\limitset_{\Gamma}\subseteq\ovl{\fixset}_{\Gamma}$, where $\ovl{\cdot}$ 
denotes the topological closure.
These sets are, in general, not equal, and the difference
$\ovl{\fixset}_{\Gamma}-\limitset_{\Gamma}$ is a discrete set,
consisting of fixed points of elliptic elements in $\Gamma$.
Now define:
$$
\tree_{\Gamma}=\tree(\limitset_{\Gamma})\quad
\textrm{and}\quad
\tree^{\ast}_{\Gamma}=\tree(\ovl{\fixset}_{\Gamma}).
$$
Clearly, we have $\tree_{\Gamma}\subseteq\tree^{\ast}_{\Gamma}$.
It is also clear that, for an inclusion $\Gamma_1\subseteq\Gamma_2$ of
finitely generated discrete subgroups, we have inclusions of trees
$\tree_{\Gamma_1}\subseteq\tree_{\Gamma_2}$ and
$\tree^{\ast}_{\Gamma_1}\subseteq\tree^{\ast}_{\Gamma_2}$.
The trees $\tree_{\Gamma}$ and $\tree^{\ast}_{\Gamma}$ admit canonically an
action by $\Gamma$ without inverstion.

\vspace{1ex}
\paragraph\label{exa-tree}\
{\bf Examples.}\ (1) If $\Gamma$ is a finite subgroup, then
$\tree_{\Gamma}$ is empty. The notion of the the other tree
$\tree^{\ast}_{\Gamma}$ fits in with the following concept: 
For an elliptic element $\gamma\in\PGL(2,K)$ with the fixed points
$z,w\in\P^1(K)$, we set
$$
M(\gamma)=\,]z,w[,
$$
and call it the {\it mirror} of $\gamma$ (this definition of mirror
slightly differs from that in \cite[(2.3)]{CKK99}; see Lemma
\ref{lem-fixedlocus} below). 
Then the tree $\tree^{\ast}_{\Gamma}$ is the minimal one which contains all
the mirrors of elements ($\neq 1$) in $\Gamma$.

(2) If $\Gamma$ is a free subgroup (i.e., so-called, {\it
Schottky group}), then the trees $\tree_{\Gamma}$ and
$\tree^{\ast}_{\Gamma}$ coincide with each other, and with the Gerritzen-van 
der Put tree $\tree^{\mathrm{GvdP}}_{\Gamma}$ (\cite[I.2.6]{GvP80}),
originally introduced by Mumford (\cite{Mum72}); indeed, in this case,
it is well-known that the set of ends of the latter tree recovers the
set of limit points (cf.\ \cite[(1.19)]{Mum72}).

(3) In general, we have
$\tree^{\mathrm{GvdP}}(\limitset_{\Gamma})=\tree_{\Gamma}$ 
(this follows easily from \cite[I.3.1\ (1)]{GvP80}), and the other tree 
$\tree^{\ast}_{\Gamma}$ is the minimal one containing $\tree_{\Gamma}$ and
all the mirrors of elliptic elements in $\Gamma$.

\vspace{1ex}
\paragraph\label{rem-mirror}\
{\bf Remark.}\ (1)\ The idea of the terminology ``mirror'' stems from the
analogy to reflection mirrors in the theory of reflection groups. In
fact, any elliptic element fixes its mirror pointwise, and ``rotates''
the other parts (cf.\ Lemma \ref{lem-fixedlocus}).

(2)\ Let $\gamma,\theta\in\Gamma$ be elliptic elements. Then the mirrors 
$M(\gamma)$ and $M(\theta)$ shares an end (i.e.,\
$M(\gamma)\bigcap M(\theta)$ contains a half-line) if and
only if $\langle\gamma,\theta\rangle$ is a cyclic group.
This follows easily from (\ref{para-recall}.3). In particular, mirrors
are in bijection with maximal finite cyclic subgroups in $\Gamma$.

\vspace{1ex}
\paragraph\label{lem-fixedlocus}\
{\bf Lemma.}\ {\sl Let $n$ be the order of $\gamma$, and set
$G=\langle\gamma\rangle$. 

(1) Let $v_0\in M(\gamma)$. If $(n,p)=1$, then $G$ acts
freely on the $q-1$ vertices adjacent to $v_0$ not lying on
$M(\gamma)$, where $q$ is the number of elements in the residue field
$k$. 

(2) Suppose $n=p^r$ for $r\geq 1$, and set
$s=\valuation(\zeta_{p^r}-1)$, where $\zeta_{p^r}$ is a primitive
$p^r$-th root of unity, and $\valuation$ is the normalized (i.e.\
$\valuation(\pi)=1$) valuation.
Then a vertex $v\in\tree_K$ is fixed by $G$
if and only if $0\leq d(v,M(\gamma))\leq s$.}

\vspace{1ex}
\pf
We may assume that $\gamma\colon z\mapsto\zeta_nz$, where $z$ is the
inhomogeneous coordinate. 
(1) follows from the fact that the adjacent vertices are in
canonical one-to-one correspondence with points in $\P^1(k)$.
(2) is due to an easy calculation collaborated with the following fact:
Let $v_0=[\O_KX_0+\O_KX_1]$ and
$v_1=[\O_K(X_0+u_0X_1)+\O_K\pi X_1]$ with $u_0\not\equiv 0\
\mathrm{mod}\ \pi$. If $v$ is a vertex such that the path
$[v_0,v]$ contains $v_1$, then
$v=[\O_K(X_0+(\sum^{d-1}_{i=0}u_i\pi^i)X_1)+\O_K\pi^dX_1]$, 
where $d=d(v,v_0)$.
\qed

\vspace{1ex}
\paragraph\label{pro-minimal}\
{\bf Proposition.}\ 
{\sl For a finitely generated discrete subgroup $\Gamma$ in $\PGL(2,K)$, 
the tree $\tree_{\Gamma}$ is minimal among the subtrees in $\tree_K$
acted on by $\Gamma$.}

\vspace{1ex}
\pf
If $\Gamma$ does not contain a hyperbolic element, then $\tree_{\Gamma}$ 
is empty and the proposition is vacuous. 
Let $\tree$ be a subtree in $\tree_K$ acted on by $\Gamma$.
It is well-known that the set of limit points $\limitset_{\Gamma}$ is
the topological closure of the set of fixed points of hyperbolic
elements. 
By Lemma \ref{lem-compact}, it suffices to show that, for any hyperbolic
element $\gamma\in\Gamma$, the tree 
$\tree$ contains the apartent connecting the fixed points of $\gamma$.
Let $v\in\Vert(\tree)$. Since $\gamma$ does not fix any vertex in
$\tree_K$, the $v_n=\gamma^nv$ for $n\in\Z$ are all distinct. For each
$n$, let $u_n$ be the vertex determined by
$[v_{n-1},v_n]\bigcap[v_n,v_{n+1}]=[v_n,u_n]$.
The vertices $u_n$ are also all distinct. 
Then $\gamma[u_{n-1},u_n]=[u_n,u_{n+1}]$, and hence $\gamma$ fixes two
ends of the apartment $\bigcup_{n\in\Z}[u_n,u_{n+1}]$ in $\tree$.
\qed

\vspace{1ex}
\paragraph\label{para-quotient}\
{\sl Quotient graphs.}\ 
We denote by $T_{\Gamma}$ and $T^{\ast}_{\Gamma}$ the quotient
graph of $\tree_{\Gamma}$ and $\tree^{\ast}_{\Gamma}$, respectively, by
$\Gamma$. 
The quotient maps of these trees are, by slight abuse of notation,
both denoted by $\varrho_{\Gamma}$.
Let $\Omega_{\Gamma}=\P^{1,\mathrm{an}}_K-\limitset_{\Gamma}$, the
corresponding analytic domain, and
$\varpi_{\Gamma}\colon\Omega_{\Gamma}\rightarrow
\Gamma\backslash\Omega_{\Gamma}$ 
the quotient map.
It is well-known that the graph $T_{\Gamma}$ is finite, and that the
analytic space 
$\Gamma\backslash\Omega_{\Gamma}$ is the analytification of a
non-singular projective curve.
Ramification points of $\varpi_{\Gamma}$ are fixed points of elliptic
elements, or equivalently, points in
$\ovl{\fixset}_{\Gamma}-\limitset_{\Gamma}$. 
This leads to the following statement (cf.\ \cite{vdP97}):

\vspace{1ex}
\paragraph\label{pro-branch}\
{\bf Proposition.}\ 
{\sl There exist canonical bijections, compatible with the quotient
maps, 
$$
\begin{array}{ccc}
\left\{
\begin{minipage}{4cm}
\setlength{\baselineskip}{.85\baselineskip}
\begin{small}
{\slshape Ramification points in $\Omega_{\Gamma}$ of the map
$\varpi_{\Gamma}$}
\end{small}
\end{minipage}
\right\}
&\longleftrightarrow&
\Ends(\tree^{\ast}_{\Gamma})-\Ends(\tree_{\Gamma})\\
\bigdownarrow&&\bigdownarrow\\
\left\{
\begin{minipage}{4cm}
\setlength{\baselineskip}{.85\baselineskip}
\begin{small}
{\slshape Branch points in $\Gamma\backslash\Omega_{\Gamma}$}
\end{small}
\end{minipage}
\right\}
&\longleftrightarrow&
\Ends(T^{\ast}_{\Gamma})
\rlap{.}
\end{array}
$$
Moreover, the decomposition group of a ramification point coincides with 
the stabilizer of the corresponding end.}
\qed

\vspace{1ex}
\paragraph\label{para-stabilizer}\
{\sl Stabilizers and tree of groups.}\ 
For $v\in\Vert(\tree_K)$ (resp.\ $\sigma\in\Edge(\tree_K)$) we
denote by $\Gamma_v$ (resp.\ $\Gamma_{\sigma}$) the stabilizer in
$\Gamma$ of $v$ (resp.\ $\sigma$ with orientation).
These are finite groups, for $\Gamma$ is discrete.
Now we assume that the quotient graph $T^{\ast}_{\Gamma}$ are trees.
Then by \cite[I.4.1, Prop.\ 17]{Ser80}, there exists a section
$\iota_{\Gamma}\colon
T^{\ast}_{\Gamma}\hookrightarrow\tree^{\ast}_{\Gamma}$ 
of the quotient map $\varrho_{\Gamma}$.
Such a section gives rise to the so-called {\it tree of groups} 
$(T^{\ast}_{\Gamma},\Gamma_{\bullet})$ (\cite[I.4.4, Def.\ 8]{Ser80}) 
by attaching the stabilizers $\Gamma_v$ (resp.\ $\Gamma_{\sigma}$) to
vertices $v\in\Vert(\iota_{\Gamma}(T^{\ast}_{\Gamma}))$
(resp.\ edges $\sigma\in\Edge(\iota_{\Gamma}(T^{\ast}_{\Gamma}))$).
By \cite[I.4.5, Thm.\ 10]{Ser80}, we see that $\Gamma$ is generated by
the finite subgroups $\Gamma_v$ for
$v\in\Vert(\iota_{\Gamma}(T^{\ast}_{\Gamma}))$, and is isomorphic to
the associated amalgam product
$$
\Gamma\stackrel{\sim}{\longrightarrow}\lim_{\longrightarrow}
(T^{\ast}_{\Gamma},\Gamma_{\bullet}).
$$
A similar isomorphy with 
$(T^{\ast}_{\Gamma},\Gamma_{\bullet})$ replaced by the finite
subtree of groups $(T_{\Gamma},\Gamma_{\bullet})$ is also true
by the same reasoning.

\vspace{1ex}
\paragraph\label{def-contraction}\
{\bf Definition}\ (cf.\ \cite[Prop.\ 1]{CKK99}). 
Let $(T,G_{\bullet})$ be an abstract tree of
groups, and $T'\subseteq T$ a subtree.
Then the induced tree of groups $(T',G_{\bullet})$ is said to be
a {\it contraction} of $(T,G_{\bullet})$ if the following conditions
are satisfied:
\begin{condlist}
\item[(1)] $\Ends(T')=\Ends(T)$.
\item[(2)] For every vertex $v$ of $T-T'$ the stabilizers of vertices on
           the path from $v$ to $v'$ are ordered increasingly with
           respect to inclusion upon approaching $T'$, where $v'$ is the
           vertex in $T'$ nearest to $v$. 
\end{condlist}

\vspace{1ex}
If $(T',G_{\bullet})$ is a contraction of $(T,G_{\bullet})$, then, by
(1), $T$ differs from $T'$ only by (possibly infinitely many) finite
pieces, and (2) means that these pieces are inessential; in
particular, the associated amalgams coincide with each other. 

\vspace{1ex}
\paragraph\label{lem-contraction}\
{\bf Lemma.}\ {\sl Let $\tree'\subseteq\tree\subseteq\tree_K$ be
inclusions of trees into the Bruhat-Tits tree $\tree_K$, and
$\Gamma\subset\PGL(2,K)$ a finitely generated discrete subgroup.
Suppose that both $\tree'$ and $\tree$ are stable under the action of
$\Gamma$. Let $T=\Gamma\backslash\tree$ and
$T'=\Gamma\backslash\tree'$, and $(T,\Gamma_{\bullet})$ and
$(T',\Gamma_{\bullet})$ the trees of groups induced by a section
$\iota\colon T\hookrightarrow\tree$. 
Suppose that the inclusion $T'\hookrightarrow T$ 
gives the bijection between the set of ends. Then $(T',\Gamma_{\bullet})$ is a
contraction of $(T,\Gamma_{\bullet})$.}

\vspace{1ex}
\pf
This is a slight generalization of \cite[Proposition 1]{CKK99}, and can
be proven by the same argument as in $\cite[(3.6)]{CKK99}$.
\qed

\vspace{2ex}
\sectioning\label{section-realization}
{\bf Realization of tree of groups}

\vspace{1ex}
\paragraph\ 
As we saw in \ref{para-stabilizer} any finitely generated discrete
subgroup $\Gamma\subset\PGL(2,K)$ such that $T^{\ast}_{\Gamma}$ is a
tree gives rise to a tree of groups
$(T^{\ast}_{\Gamma},\Gamma_{\bullet})$, which recovers the abstract
group isomorphic to $\Gamma$ as, so to speak, the ``fundamental group'' of  
the data $(T^{\ast}_{\Gamma},\Gamma_{\bullet})$. 
Moreover the data
$(T^{\ast}_{\Gamma},\Gamma_{\bullet})$ also recovers
$\tree^{\ast}_{\Gamma}$ as an abstract tree (cf.\ \cite[I.4.5, Thm.\
10]{Ser80}), which one can call the ``universal covering'' of
$(T^{\ast}_{\Gamma},\Gamma_{\bullet})$. 
Now the natural question rises: Given an abstract tree of groups, when
can one realize its fundamental group as a
discrete subgroup in $\PGL(2,K)$ and the universal covering 
as a subtree in $\tree_K$?
The objective of this section is to answer this question.

\vspace{1ex}
\paragraph\label{para-situation}\
Let $(T,G_{\bullet})$ be an abstract tree of groups, that is, an
abstract tree $T$ to which finite groups $G_v$ and $G_{\sigma}$ for
$v\in\Vert(T)$ and $\sigma\in\Edge(T)$ are attached; among these groups
are injective homomorphisms $G_{\sigma}\hookrightarrow G_v$ for each
pair $(v,\sigma)$ with $v\vdash\sigma$. Suppose that we are given
embeddings $G_v\hookrightarrow\PGL(2,K)$ for any $v\in\Vert(T)$
compatible with each $G_{\sigma}\hookrightarrow G_v$ for 
$v\vdash\sigma$. Such embeddings, provided that $K$ is large enough,
gives rise to subtrees $\tree^{\ast}_{G_v}$ as in
\ref{para-discretegroup}.
Set 
$$
\til{\tree}_{G_{\bullet}}=\textrm{the minimal subtree in $\tree_K$
containing all $\tree^{\ast}_{G_v}$}.
$$
The set of ends in $\til{\tree}_{G_{\bullet}}$ is, therefore, the union
of the set of ends in $\tree^{\ast}_{G_v}$ for $v\in\Vert(T)$.
This tree is labelled by groups (not necessarily finite)
$\til{G}_{\bullet}$ as follows: For a vertex
$v\in\Vert(\til{\tree}_{G_{\bullet}})$ the group $\til{G}_v$ is the
subgroup in $\PGL(2,K)$ generated by $(G_u)_v$ (the stabilizer at $v$ by 
the action of $G_u$ on $\tree_K$) for all $u\in\Vert(T)$; the definition 
of the group $\til{G}_{\sigma}$ for $\sigma\in\Edge(\til{\tree}_{G_{\bullet}})$
is similar, which is just the intersection of $\til{G}_v$'s at the two
extremities.

\vspace{1ex}
\paragraph\label{def-admissible}\
{\bf Definition.}\ An {\it admissible embedding} of an abstract tree of
groups $(T,G_{\bullet})$ is an embedding $\iota\colon T\hookrightarrow\tree_K$ 
of trees together with embeddings $G_v\hookrightarrow\PGL(2,K)$ for any
$v\in\Vert(T)$ compatible with each $G_{\sigma}\hookrightarrow G_v$ for
any $v\vdash\sigma$ such that the following conditions are satisfied:
\begin{condlist}
\item[(1)] $\iota(T)\subset\til{\tree}_{G_{\bullet}}$.
\item[(2)] For any $v\in\Vert(T)$ and $\gamma\in G_v$ ($\gamma\neq 1$),
           there exists $\delta\in\Gamma$ such that
           $M(\delta\gamma\delta^{-1})\bigcap\iota(T)$ contains an edge, 
           where $\Gamma$ is the subgroup in $\PGL(2,K)$ generated by
           all $G_v$ for $v\in\Vert(T)$.
\item[(3)] $\til{G}_{\iota(v)}=G_v$ for any $v\in\Vert(T)$. 
\item[(4)] $\til{G}_{\iota(\sigma)}=G_{\sigma}$ for any
           $\sigma\in\Edge(T)$.
\item[(5)] For any $v\in\Vert(T)$, we have $\Star_v(T)\cong G_v\backslash
           (G_v\cdot\Star_{\iota(v)}(\til{\tree}_{G_{\bullet}}))$ by the
           composite of $\iota$ followed by the projection.
\end{condlist}
The last condition means that $T$ behaves locally like a fundamental
domain at each vertex.

\vspace{1ex}
\paragraph\label{def-staradmissible}\
{\bf Definition.}\ An abstract tree of groups $(T,G_{\bullet})$ is said
to be {\it $\ast$-admissible} if it has an admissible embedding and the
associated amalgam $\lim_{\rightarrow}(T,G_{\bullet})$ is finitely
generated.

\vspace{1ex}
\paragraph\label{lem-admissible}\
{\bf Lemma.}\ 
{\sl If $\Gamma\in\PGL(2,K)$ is a finitely generated discrete
subgroup such that $T^{\ast}_{\Gamma}$ is a tree, then
$(T^{\ast}_{\Gamma},\Gamma_{\bullet})$ by a section
$\iota_{\Gamma}\colon
T^{\ast}_{\Gamma}\hookrightarrow\tree^{\ast}_{\Gamma}$ 
is $\ast$-admissible.}

\vspace{1ex}
\pf
By \ref{para-stabilizer},
$\lim_{\rightarrow}(T,G_{\bullet})\cong\Gamma$, and is finitely
generated.
Clearly, we have
$\iota_{\Gamma}(T^{\ast}_{\Gamma})\subset\til{\tree}_{\Gamma_{\bullet}}$.
(\ref{def-admissible}.3) and
(\ref{def-admissible}.4) are obvious.
Let $v\in\Vert(T^{\ast}_{\Gamma})$. Then
$\Star_v(T^{\ast}_{\Gamma})\cong
\Gamma_v\backslash\Star_v(\tree^{\ast}_{\Gamma})$ obviously holds. But
since $\iota_{\Gamma}(T^{\ast}_{\Gamma})
\subset\til{\tree}_{\Gamma_{\bullet}}\subset\tree^{\ast}_{\Gamma}$, we 
have (\ref{def-admissible}.5). 
Finally, for $\gamma\in G_v$ with $\gamma\neq 1$, since $T^{\ast}$ is a
fundamental domain in $\tree^{\ast}_{\Gamma}$, there exists
$\delta\in\Gamma$ such that $M(\delta\gamma\delta^{-1})\bigcap T$ is
non-empty, containing a vertex $w$. Due to (\ref{def-admissible}.5),
of which we have proved the validity, one can further make a twist by
$\chi\in G_w$ so that $M(\chi\delta\gamma\delta^{-1}\chi^{-1})\bigcap T$ 
contains an edge.
\qed

\vspace{1ex}
Note that $(T_{\Gamma},\Gamma_{\bullet})$ is not $\ast$-admissible,
since it does not satisfy (\ref{def-admissible}.5). What we are to show
is that the converse of the above lemma in a certain sense:

\vspace{1ex}
\paragraph\label{thm-realization}\
{\bf Theorem.}\ 
{\sl Let $(T,G_{\bullet})$ be a $\ast$-admissible tree of groups and 
$\iota\colon T\hookrightarrow\til{\tree}_{G_{\bullet}}$ with 
$\{G_v\hookrightarrow\PGL(2,K)\}_{v\in\Vert(T)}$ an admissible
embedding. Let $\Gamma$ be the subgroup in $\PGL(2,K)$ generated by all
$G_v$ for $v\in\Vert(T)$ and set
$$
\tree^{\ast}=\bigcup_{\gamma\in\Gamma}\gamma\cdot\iota(T)
$$
in $\tree_K$. Then:

\vspace{1ex}
(1) The group $\Gamma$ is a finitely generated discrete subgroup in
$\PGL(2,K)$ isomorphic to $\lim_{\rightarrow}(T,G_{\bullet})$.

\vspace{1ex}
(2) The subset $\tree^{\ast}$ in $\tree_K$ is a tree and
$\Gamma\backslash\tree^{\ast}\cong T$.

\vspace{1ex}
(3) The embedding $\iota$ gives a section $T\hookrightarrow\tree^{\ast}$ 
by which the induced tree of groups $(T,\Gamma_{\bullet})$ equals to 
$(T,G_{\bullet})$.

\vspace{1ex}
Moreover, if $\tree^{\ast}_{\Gamma}$ is the tree associated to $\Gamma$ as
in \ref{para-discretegroup}, then
$\tree^{\ast}_{\Gamma}\subseteq\tree^{\ast}$, and the induced inclusion
$T^{\ast}_{\Gamma}\hookrightarrow T$ enjoys the following:

\vspace{1ex}
(4) The induced inclusion
$\Ends(T^{\ast}_{\Gamma})\hookrightarrow\Ends(T)$ is a bijection.

\vspace{1ex}
(5) The section $\iota$ restricts to a section
$T^{\ast}_{\Gamma}\hookrightarrow\tree^{\ast}_{\Gamma}$ by which the
induced tree of groups $(T^{\ast}_{\Gamma},\Gamma_{\bullet})$ is the
restriction of $(T,\Gamma_{\bullet})=(T,G_{\bullet})$.

\vspace{1ex}
(6) The tree of groups $(T^{\ast}_{\Gamma},\Gamma_{\bullet})$ is a
contraction of $(T,\Gamma_{\bullet})=(T,G_{\bullet})$.}

\vspace{1ex}
\paragraph\label{para-length}\
To prove the theorem, we need several lemmas.
In the sequel, we regard $T$ as a subtree in $\tree_K$ by $\iota$; also, 
we can simply write $G_v$ and $G_{\sigma}$ instead of $\til{G}_v$ and
$\til{G}_{\sigma}$, respectively, because of (\ref{def-admissible}.3)
and (\ref{def-admissible}.4).
Trees are often regarded as metric spaces by geometric realization
(cf.\ \ref{para-treenotation}).

Let $\Gamma$ be as in the theorem. Any element $\gamma\in\Gamma$ is
expressed as 
$\gamma=\alpha_1\cdots\alpha_m$ 
with $\alpha_i\in G_{v_i}$ for some $v_i\in\Vert(T)$ for
$i=1,\ldots,m$. 
The length of $\gamma$ is the minimal $m$ in among all such expressions
as above. Let $\Gamma^{(m)}$ be the set of all elements in $\Gamma$ of
length $m$. Obviously, $\Gamma^{(0)}=\{1\}$ and
$\Gamma^{(1)}=\bigcup_{v\in\Vert(T)}G_v-\{1\}$.

\vspace{1ex}
\paragraph\label{lem-nontrivial}\
{\bf Lemma.}\ {\sl There is no half-line in $T$ on which only trivial
groups are attached to vertices and edges.}

\vspace{1ex}
\pf
From the definition of the tree $\til{\tree}_{G_{\bullet}}$, it follows
that $\Ends(\til{\tree}_{G_{\bullet}})$ is, regarded as a subset in
$\P^1(K)$, the set of fixed points of elements in $\Gamma^{(1)}$.
Then the lemma follows from (\ref{def-admissible}.1),
(\ref{def-admissible}.3), and (\ref{def-admissible}.4).
(Note that if $G_u=1$ then $\tree^{\ast}_{G_u}=\emptyset$.)
\qed

\vspace{1ex}
\paragraph\label{lem-arcwise}\
{\bf Lemma.}\ {\sl The subset $\tree^{\ast}$ in $\tree_K$ is arcwise 
connected, i.e.\ a subtree.}

\vspace{1ex}
\pf
For any $\gamma=\alpha_1\cdots\alpha_m\in\Gamma^{(m)}$, set
$\gamma_i=\alpha_1\cdots\alpha_i$ for $i=1,\ldots,m$ (set $\gamma_0=1$).
Then
$\gamma_iT\bigcap\gamma_{i+1}T=\gamma_i(T\bigcap\alpha_{i+1}T)
\neq\emptyset$ for $i=0,\ldots,m-1$.
Hence a point in $\gamma T$ can be connected by a path with a point in
$T$.
\qed

\vspace{1ex}
\paragraph\label{lem-fundamental}\
{\bf Lemma.}\ 
{\sl For $\gamma\in\Gamma^{(1)}$ and $v\in\Vert(T)$, $\gamma
v\in\Vert(T)$ implies $\gamma v=v$.} 

\vspace{1ex}
\pf
Take $u\in\Vert(T)$ such that $\gamma\in G_u$. Let $w\in\Vert(T)$ be the 
vertex determined by 
$[u,v]\bigcap[u,\gamma v]=[u,w]$. Since $\gamma u=u$ we have $\gamma
w=w$, i.e., $\gamma\in G_w$. If $\gamma v\neq v$, then the segments
$[w,v]$ and $[w,\gamma v]$ in $T$
contain edges, different from each other, emanating from $w$ which are
in the same $\gamma$-orbit. But this contradicts
(\ref{def-admissible}.5).  
\qed

\vspace{1ex}
\paragraph\label{lem-disjoint}\
{\bf Lemma.}\ {\sl Let $\gamma\in\Gamma^{(m)}$ with $m>1$.
Then $T\bigcap\gamma T=\emptyset$. 
Moreover, if $\gamma=\alpha_1\cdots\alpha_m$ a minimal expression, and 
$v_1,v_2\in\Vert(T)$ with $\alpha_1\in G_{v_1}$ and $\alpha_2\in
G_{v_2}$ are chosen so that $d(v_1,v_2)$ is minimal, then the geodesic
path connecting $T$ and $\gamma T$ contains $\alpha_1 v_2$.}

\vspace{1ex}
\pf
The proof is done by induction with respect to $m$.
First we show the lemma in $m=2$; $\gamma=\alpha_1\alpha_2$.
Since $\alpha_1\not\in G_{v_2}$, $\alpha_1 v_2=\gamma
v_2\not\in\Vert(T)$ (due to Lemma \ref{lem-fundamental}).
Now suppose $v\in T\bigcap\gamma T$. Then $[v,\alpha_1 v_2]=[v,\gamma
v_2]\subseteq\gamma T$. Due to the minimality of $d(v_1,v_2)$, we have
$[v_1,\alpha_1 v_2]\bigcap T=\{v_1\}$. Hence the geodesic path
connecting $\gamma v_2$ with the vertex $v$ in $T$ contains $v_1$; in
particular, $v_1\in T\bigcap\gamma T$.
This means $\gamma^{-1}v_1\in\Vert(T)$, while $\gamma^{-1}v_1\neq v_1$
(since $\gamma\not\in\Gamma^{(1)}$).  
But $\gamma^{-1}v_1=\alpha_2^{-1}v_1$ leads to contradiction to Lemma
\ref{lem-fundamental}.  
Therefore, $T\bigcap\gamma T=\emptyset$.

Due to the minimality of $d(v_1,v_2)$, $[v_1,\alpha_1 v_2]\bigcap
T=\{v_1\}$ and $[v_2,\alpha_2^{-1}v_1]\bigcap T=\{v_2\}$. This last
equality gives $[\alpha_1 v_2,v_1]\bigcap\gamma T=\{\alpha_1 v_2\}$.
Hence the segment $[v_1,\alpha_1 v_2]$ is the geodesic path connecting
$T$ and $\gamma T$, which contains $\alpha_1 v_2$.

For $m>2$, we set $\gamma'=\alpha_1\cdots\alpha_{m-1}$.
Take $v_m\in\Vert(T)$ such that $\alpha_m\in G_{v_m}$.
Then $\gamma v_m=\gamma' v_m\in\gamma' T\bigcap\gamma T$.
Suppose $v\in T\bigcap\gamma T$.
By induction, $[\gamma' v_m,v]$ contains $\alpha_1 v_2$.
But this segment $[\gamma' v_m,v]$ is included in $\gamma T$, which
means $\gamma T\bigcap\alpha_1 T\neq\emptyset$. This contradics to the
induction hypothesis, since $\gamma T\bigcap\alpha_1
T=\alpha_1(\alpha_2\cdots\alpha_mT\bigcap T)$.
Hence $T\bigcap\gamma T=\emptyset$.
Since the geodesic connecting $\gamma T$ and $T$ contains that
connecting $\gamma' T$ and $T$ (since $\gamma' T\bigcap\gamma
T\neq\emptyset$), in particular, it contains $\alpha_1 v_2$.
\qed

\vspace{1ex}
\paragraph\label{cor-disjoint1}\
{\bf Corollary.}\ {\sl For $\gamma\in\Gamma$ ($\gamma\neq 1$),
$T\bigcap\gamma T\neq\emptyset$ if and only if $\gamma\in\Gamma^{(1)}$.}
\qed

\vspace{1ex}
\paragraph\label{cor-disjoint1.5}\
{\bf Corollary.}\ {\sl For $\gamma\in\Gamma$ and $v\in\Vert(T)$,
$\gamma v\in\Vert(T)$ imples $\gamma\in G_v$.}
\qed

\vspace{1ex}
\paragraph\label{cor-disjoint2}\
{\bf Corollary.}\ {\sl For $\gamma,\delta\in\Gamma$, and suppose
$T\bigcap\gamma T\bigcap\gamma\delta T\neq\emptyset$. Then, for any
$v\in T\bigcap\gamma T\bigcap\gamma\delta T$, we have $\gamma,\delta\in
G_v$.} 

\vspace{1ex}
\pf
By Lemma \ref{lem-disjoint}, $\gamma$, $\delta$, and $\gamma\delta$
are in $\Gamma^{(1)}\bigcup\{1\}$.
Since $\gamma^{-1}v\in\Vert(T)$, $\gamma^{-1}v=v$ (by Lemma
\ref{lem-fundamental}), which gives $\gamma\in G_v$. 
Similarly, we get $\delta\in G_v$.
\qed

\vspace{1ex}
By these corollaries and \cite[Appendix, pp.\ 30--31]{Ser80}, we have:

\vspace{1ex}
\paragraph\label{cor-disjoint3}\
{\bf Corollary.}\ {\sl Let $F$ be the free group with basis $X_{\alpha}$ 
indexed by $\alpha\in\Gamma^{(1)}\bigcup\{1\}$, and $\varphi\colon
F\rightarrow\Gamma$ the natural homomorphism $X_{\alpha}\mapsto\alpha$.
Then $\mathrm{Ker}\varphi$ is the normal subgroup generated by
$X_{\alpha}X_{\beta}(X_{\alpha\beta})^{-1}$ for all $(\alpha,\beta)$
such that $\alpha,\beta\in G_v$ for some $v\in\Vert(T)$.}
\qed

\vspace{1ex}
\paragraph\label{cor-disjoint4}\
{\bf Corollary.}\ {\sl The natural homomorphism 
$$
\lim_{\longrightarrow}(T,G_{\bullet})\longrightarrow\Gamma
$$
is an isomorphism.}

\vspace{1ex}
\pf
It suffices to show that the kernel of the homomorphism
$F\rightarrow\lim_{\longrightarrow}(T,G_{\bullet})$ sending
$X_{\alpha}\mapsto\alpha$ is the normal subgroup generated by 
$X_{\alpha}X_{\beta}(X_{\alpha\beta})^{-1}$ for $\alpha,\beta\in
G_v\subset\lim_{\longrightarrow}(T,G_{\bullet})$ with some $v\in\Vert(T)$.
But this is obvious from the definition of amalgams.
\qed

\vspace{1ex}
\paragraph\label{para-proof}\
{\sl Proof of Theorem \ref{thm-realization}.}\ 
Let $X$ be the abstract tree (``universal covering'' of
$(T,G_{\bullet})$) as in \cite[I.4.5\ Theorem 9]{Ser80}. We first show 
that our tree $\tree^{\ast}$ and $X$ are $\Gamma$-equivariantly isomorphic. 
To see this, it suffices to show that 
$$
\Vert(\tree^{\ast})\ (=\Gamma\Vert(T))\cong
\coprod_{v\in\Vert(T)}G_T/G_v,
$$
where $G_T=\lim_{\longrightarrow}(T,G_{\bullet})$, and that the similar
equality holds also for the set of oriented edges. 
But these follow from Corollary \ref{cor-disjoint4} and Corollary
\ref{cor-disjoint1.5}.

Then it follows from \cite[I.4.5\ Theorem 9]{Ser80} that $T$ is a
fundamental domain for $\tree^{\ast}$ modulo $\Gamma$, and the
stabilizers $\Gamma_v$ ($v\in\Vert(T)$) and $\Gamma_{\sigma}$
($\sigma\in\Edge(T)$) are equal to $G_v$ and $G_{\sigma}$, repectively.
In particular, by Lemma \ref{lem-discrete}, $\Gamma$ is discrete in
$\PGL(2,K)$. Therefore, (1) and (2) have been proved.
The embedding $\iota$ obviously gives a section
$T\rightarrow\tree^{\ast}$, and hence, we have (3).

We proceed to the proof of (4)$\sim$(6).
First we are going to show $\tree^{\ast}_{\Gamma}\subseteq\tree^{\ast}$.
In view of Proposition \ref{pro-minimal}, it suffices to show that the
mirror of any elliptic element $\gamma\in\Gamma$ is contained in
$\tree^{\ast}$.
Since $\gamma$ is $\Gamma$-conjugate to an element in $\Gamma^{(1)}$
(\cite[I.1.3\ Corollary 1]{Ser80}), by 
(\ref{def-admissible}.2), we may assume that $\gamma\in G_v$ for a
vertex $v\in\Vert(T)$ and $M(\gamma)\bigcap T$ contains an edge.
If $M(\gamma)\subseteq T$, there is nothing to prove. 
Otherwise, $M(\gamma)\bigcap T$ is either a 
half-line, a segment of finite length.
If it is a half-line $\ell=[u,\varepsilon[$, then let us denote the
other ``half'' by $\ovl{\ell}=[u,\ovl{\varepsilon}[$
($M(\gamma)=\ell\bigcup\ovl{\ell}$). Since $\varepsilon$ and 
$\ovl{\varepsilon}$ are in the same orbit by the action of $G_u$ on
$\P^1(K)$, we have $\delta\in G_u$ such that 
$\delta\ell=\ovl{\ell}$. Hence $M(\gamma)\subseteq\tree^{\ast}$.
Suppose $M(\gamma)\bigcap T$ is a segment $[u,w]$.
Let $[w_1,u]$ be in $M(\gamma)$ such that $d(w_1,u)=d(u,w)$ and 
$[w_1,u]\bigcap[u,w]=\{u\}$.
Then in $\tree^{\ast}_{G_u}$ these two segements are in the same
$G_u$-orbit by the same reasoning as above for two half-lines in
$M(\gamma)$ starting at $u$ extending these two segments.
We can find $\delta_1\in G_u$ such that
$\delta_1[u,w]=[w_1,u]$. 
We can do the same for $[w_1,w]$ looking at $w_1=\delta_1(w)$ and
$\delta_1\tree^{\ast}_{G_w}=\tree^{\ast}_{G_{w_1}}$.
We find $\delta_2\in G_{w_1}\subset\Gamma$ such that
$\delta_2[w_1,u]=[w_2,w_1]$ in $M(\gamma)$. Repeating this, we can
inductively find $w_n$ such that $[w_n,w]$ is in $M(\gamma)$ and
$\delta_n\in\Gamma$ such that $\delta_n[w_n,w]=[w_{n+1},w_n]$.
These segments are in $\tree^{\ast}$, and the union of them is a
half-line starting at $w$ contained in $M(\gamma)$. Similarly, we can
find the other half in $\tree^{\ast}$. 
Hence we have shown that all the mirrors of elliptic elements in
$\Gamma$ appear in $\tree^{\ast}$, thereby
$\tree^{\ast}_{\Gamma}\subseteq\tree^{\ast}$. 

Next we claim
that $\Ends(\tree^{\ast})$ is, as a subset of $\P^1(K)$, equal to
$\ovl{\fixset}_{\Gamma}$ (cf.\ \ref{para-discretegroup}).
It follows from $\tree^{\ast}_{\Gamma}\subseteq\tree^{\ast}$ that 
$\ovl{\fixset}_{\Gamma}$ is contained in $\Ends(\tree^{\ast})$. 
If there exists a half-line $\ell$ in $\tree^{\ast}$ pointing to
$\varepsilon\not\in\ovl{\fixset}_{\Gamma}$, then, replaced by a
subhalf-line if necessary, $\ell$ contains no vertex with a non-trivial
stabilizer, and hence is mapped to a half-line in $T$ on 
which the stabilizers of vertices and edges are all trivial groups; but
this contradicts Lemma \ref{lem-nontrivial}. 
Hence we get
$\Ends(\tree^{\ast})=\ovl{\fixset}_{\Gamma}=\Ends(\tree^{\ast}_{\Gamma})$, 
and we obtain the bijection in (4) by taking quotient by
$\Gamma$. 

Since the diagram of morphism of trees
$$
\begin{array}{ccc}
\tree^{\ast}_{\Gamma}&\longhookrightarrow&\tree^{\ast}\\
\bigdownarrow&&\bigdownarrow\\
T^{\ast}_{\Gamma}&\longhookrightarrow&T
\end{array}
$$
is cartesian (e.g.\ in the category of metric spaces), it follows that
the section $\iota$ restricts to a section of $T^{\ast}_{\Gamma}$ into 
$\tree^{\ast}_{\Gamma}$. 
The other part of (5) is clear. 
(6) is due to Lemma \ref{lem-contraction}.
\qed

\vspace{1ex}
To conclude this section, we insert herein a corollary to Theorem
\ref{thm-realization} useful for application. 
Let $(T,G_{\bullet})$ be a $\ast$-admissible tree of groups.
In view of Proposition \ref{pro-branch} the ends of $T$ are in bijection 
with branch points of
$\Omega_{\Gamma}\rightarrow\Gamma\backslash\Omega_{\Gamma}\cong\P^1_K$. 
Hence $\Ends(T)$ is a finite set.
Let $\varepsilon\in\Ends(T)$. Since $\lim_{\longrightarrow}(T,G_{\bullet})$
is finitely generated, we can find a half-line in $T$ converging to
$\varepsilon$ such that the attached groups are ordered decreasingly
with respect to inclusion upon approaching $\varepsilon$. We denote by
$G_{\varepsilon}$ the intersection of these groups, and call the
stabilizer of $\varepsilon$. This is not a trivial group due to
Lemma \ref{lem-nontrivial}.

\vspace{1ex}
\paragraph\label{cor-realization}\
{\bf Corollary.}\ {\sl Let $(T,G_{\bullet})$ be a $\ast$-admissible tree
of groups. Then, for $\varepsilon\in\Ends(T)$, the group
$G_{\varepsilon}$ is a finite cyclic group. Let
$\Ends(T)=\{\varepsilon_1,\ldots,\varepsilon_n\}$, and $o(i)$ the order
of $G_{\varepsilon_i}$ for $i=1,\ldots,n$.
Then there exists a finitely generated discrete subgroup $\Gamma$ in
$\PGL(2,K)$ isomorphic to $\lim_{\longrightarrow}(T,G_{\bullet})$ such
that $\Gamma\backslash\Omega_{\Gamma}\cong\P^1_K$ and the quotient map
$\varpi_{\Gamma}\colon\Omega_{\Gamma}\rightarrow
\Gamma\backslash\Omega_{\Gamma}$ branches over $n$ points with branching 
degrees $o(\varepsilon_1),\ldots,o(\varepsilon_n)$.}

\vspace{1ex}
\pf
All these are clear by the theorem and Proposition \ref{pro-branch}.
\qed

\vspace{2ex}
\sectioning\label{section-application}
{\bf Examples}

\vspace{1ex}
\paragraph\label{para-free}\
{\sl Free product} (cf.\ \cite[\S 11]{Her78}).\ 
Let us begin with a simple example. Let $(T,G_{\bullet})$
be the tree of groups as drawn in Figure 1.
$$
\begin{picture}(150,150)(0,-10)
\put(35,75){\vector(0,1){50}}
\put(35,75){\vector(0,-1){50}}
\put(115,75){\vector(0,1){50}}
\put(115,75){\vector(0,-1){50}}
\put(35,75){\line(1,0){80}}
\put(35,75){\circle*{4}}
\put(115,75){\circle*{4}}
\put(20,72){$\scriptstyle{Z_n}$}
\put(120,72){$\scriptstyle{Z_m}$}
\put(74,78){$\scriptstyle{1}$}
\put(54,0){{\small{\sc Figure 1}}}
\end{picture}
$$

\vspace{-2ex}\noindent
Here the four arrows stand for ends; to the left (resp.\ right) vertical
line only the cyclic group $Z_n$ of order $n$ (resp.\ $Z_m$ of order
$m$) is attached, while the groups attached to the vertices on the
horizontal line, except for its extremities, are all trivial groups.

Let us show that, provided that $n$ and $m$ are prime to $p$ and that
the length of the horizontal line is even, the $(T,G_{\bullet})$ is
$\ast$-admissible (the assumption on the length of the horizontal line
is by no means essential, since one can always attain it by replacing
$K$ by a ramified quadratic extension): 
Let $K$ be a finite extension of $\Q_p$ containing $\zeta_n$ (resp.\
$\zeta_m$), a primitive $n$-th (resp.\ $m$-th) root of unity, 
$\pi\in\O_K$ a prime element, and $r$ the half of the length of
the horizontal line.
The embedding $\iota\colon T\rightarrow\tree_K$ is the one determined as 
follows: The left (resp.\ right) vertical line is mapped to the
apartment $]0,\pi^r[$ (resp.\ $]\pi^{-r},\infty[$). These apartments are 
disjoint and of distance $2r$; in fact, the segment connecting
$v_0=[\O_K\pi^r\oplus\O_K]$ and $v_1=[\O_K\oplus\O_K\pi^r]$ is the
geodesic path, which is the image of the horizontal segment in $T$.
Let $\gamma,\delta\in\PGL(2,K)$ be defined by the fractional linear
transformations:
$$
\gamma(z)=\frac{\zeta_n\pi^rz}{(\zeta_n-1)z+\pi^r},
\qquad
\delta(z)=\zeta_mz-(\zeta_m-1)\pi^{-r},
$$
where $z$ is the inhomogenious coordinate.
The element $\gamma$ (resp.\ $\delta$) is of order $n$ (resp.\ $m$), and 
has $0$ and $\pi^r$ (resp.\ $\pi^{-r}$ and $\infty$) as its fixed
points.  
Let us embedd $Z_n$ and $Z_m$ in $(T,G_{\bullet})$ by fixing
$Z_n\cong\langle\gamma\rangle$ and $Z_m\cong\langle\delta\rangle$. 
Then by Lemma \ref{lem-fixedlocus}, if $n$ and $m$ are prime to $p$, we
see $\til{\tree}_{G_{\bullet}}=\iota(T)$, and we can easily check that 
the above embedding is admissible.
Hence by Theorem \ref{thm-realization}, the subgroup
$\Gamma=\langle\gamma,\delta\rangle$ is discrete and isomorphic to the 
free product $Z_n\ast Z_m$. By Corollary \ref{cor-realization}, this
$\Gamma$ gives $\Omega_{\Gamma}\rightarrow\P^1_K$ brached over $4$
points with branching degrees $(n,n,m,m)$.

If either $n$ or $m$ are not prime to $p$, then one can modify the groups on
the horizontal line according to Lemma \ref{lem-fixedlocus} (hence $r$
should be large enough) to make it $\ast$-admissible. Also in this case
the associated group $\Gamma$ is, provided $r$ large enough, isomorphic
to the free product $Z_n\ast Z_m$.

\vspace{1ex}
\paragraph\label{para-triangle}\
{\sl Triangle group.}\ 
In this paragraph we assume $p=2$.
Let $n$ be a positive odd number, and $K$ a finite extension of $\Q_2$
containing $\zeta_n$, and $\pi\in\O_K$ a prime. 
Let $e$ be the ramification degree of $K$ over $Q_2$. 
First consider the dihedral subgroup $D_n$ generated by
$\gamma,\chi\in\PGL(2,K)$ with
$$
\gamma(z)=\zeta_nz,\qquad\chi(z)=1/z.
$$
The fixed points of $\gamma$ are $0$ and $\infty$, while those of $\chi$
are $1$ and $-1$.
One sees easily that $M(\gamma)=]0,\infty[$ and $M(\chi)=]1,-1[$
are disjoint with distance $e$.
Let $v_0$ (resp.\ $v_1$) be the vertex in $\tree_K$ which is the
similarity class of the standard lattice $\O_Ke_0+\O_Ke_1$ (resp.\ the
lattice $\O_K(e_0+e_1)+\O_K2e_1$), where $e_0=(1,0)$ and $e_1=(0,1)$.
Then the segment $[v_0,v_1]$ is the geodesic path connecting $M(\gamma)$ 
and $M(\chi)$.
A fundamental domain $T^{\ast}_{D_n}$ for $\tree^{\ast}_{D_n}$ modulo
$D_n$ is given by the union of (i) the half-line $[v_0,\infty[$, (ii) 
$M(\chi)$, and (iii) the segment $[v_0,v_1]$ (see Figure 2).
$$
\begin{picture}(150,150)(0,-10)
\put(75,55){\vector(2,-1){64}}
\put(75,55){\vector(-2,-1){64}}
\put(75,55){\circle*{4}}
\put(75,73){\circle*{4}}
\multiput(75,56)(0,4){4}{\circle*{2}}
\put(75,80){\vector(0,1){47}}
\put(75,80){\vector(0,1){3}}
\put(71,43){$\scriptstyle{Z_2}$}
\put(80,70){$\scriptstyle{D_n}$}
\put(63,72){$\scriptstyle{v_0}$}
\put(70,132){$\scriptstyle{Z_n}$}
\put(63,58){$\scriptstyle{v_1}$}
\put(1,17){$\scriptstyle{Z_2}$}
\put(143,17){$\scriptstyle{Z_2}$}
\put(28,0){{\small{\sc Figure 2}: $n$: odd}}
\end{picture}
$$

\vspace{-2ex}\noindent
The groups attached to $T^{\ast}_{D_n}$ are as follows: On $M(\chi)$ all 
vertices and edge are labelled by $Z_2$. Vertices and edges on
$[v_0,\infty[$, except for $v_0$ are labelled by $Z_n$, while $v_0$ is
by the whole $D_n$.
By Lemma \ref{lem-fixedlocus}, to vertices and edges in $[v_0,v_1]$
(denoted by the dotted segment),
except for $v_0$, the group $Z_2$ is attached. (Needless to say, they
are subgroups of $D_n$).

Now consider $\theta\in\PGL(2,K)$ elliptic of order $2m$ with $m$ odd
such that $\theta^m=\chi$, i.e., $\theta$ has the same fixed points as
$\chi$. Then we can consider the (abstract) tree of groups as in Figure
2 with $M(\chi)$ 
replaced by $M(\theta)$, which amounts to replace all the $Z_2$'s on the
lower straight-line in Figure 2 by $Z_{2m}$.
It can be checked that the resulting tree of groups $(T,G_{\bullet})$ is
$\ast$-admissible by the obvious embeddings; that $m$ is assumed to be
odd guarantees 
(\ref{def-admissible}.3) and (\ref{def-admissible}.4) (on vertices and
edges in $[v_0,v_1]$) due to Lemma \ref{lem-fixedlocus}.
Although in this case the tree $\til{\tree}_{G_{\bullet}}$ is bigger
than $T$, validity of (\ref{def-admissible}.5) follows from an argument
similar to that in the previous example, and the fact that $T$ came from 
the fundamental domain of $D_4$.

It follows therefore that for any odd numbers $n$ and $m$ there exists a 
discrete subgroup $\Gamma$ in $\PGL(2,K)$ (with $K$ sufficiently large)
isomorphic to $D_n\ast_{Z_2}Z_m$ such that the associated quotient map
$\varpi\colon\Omega_{\Gamma}\rightarrow\Gamma\backslash\Omega_{\Gamma}
\cong\P^1_K$ braches exactly above three points with branching degree 
$(n,2m,2m)$.

The assumption that $m$ is odd is actually not essential; even if $m$ is
an even number, one can modify the groups attached to $[v_0,v_1]$,
accorting to Lemma \ref{lem-fixedlocus} so
that the resulting tree of groups is $\ast$-admissible. 

\vspace{1ex}
\paragraph\label{rem-triangle}\
{\bf Remark.}\ (1)\ Note that, unless $p=2$, the above construction does 
not work any more, since the fundamental domain of $T^{\ast}_{D_n}$
looks different; more precisely, if $p>2$ then $M(\gamma)$ and
$M(\chi)$ has non-empty intersection. Hence one cannot perform the
replacement of groups as above.

(2)\ The resulting discrete group is a $p$-adic analogue of the
Schwarzian triangle groups (cf.\ \cite[\S 9]{And98}).
Our example gives an affirmative answer to Yves Andr\'e's expectation
(cf.\ \cite[9.4]{And98}) 
that there will be infinitely many non-arithmetic $p$-adic triangle
groups.

(3)\ Using the method as above, we can construct more triangle groups,
not only in $p=2$; by this, in particular, one can show that, for
$p=2,3,5$ there are infinitely many non-arithmetic triangle
groups. Actually, one can also show that for $p>5$ there is no triangle
group constructed by this method. 
In order to discuss these, as the construction in \ref{para-triangle}
indicates, one has to describe the fundamental domains for finite groups 
of other types, i.e., $D_n$ with $n$ even, tetrahedral group,
octahedral group, and icosahedral group. 
This will be done in \cite{Kat00}.

\hyphenation{Schottky-kurve}
\begin{small}

{\sc Graduate School of Mathematics, Kyushu University, Hakozaki
Higashi-ku, Fukuoka 812-8581, Japan.}
\end{small}
\end{document}